\PassOptionsToPackage{unicode}{hyperref}
\PassOptionsToPackage{hyphens}{url}
\documentclass[
]{article}
\usepackage{lmodern}
\usepackage{amssymb,amsmath}
\usepackage{ifxetex,ifluatex}
\ifnum 0\ifxetex 1\fi\ifluatex 1\fi=0 
  \usepackage[T1]{fontenc}
  \usepackage[utf8]{inputenc}
  \usepackage{textcomp} 
\else 
  \usepackage{unicode-math}
  \defaultfontfeatures{Scale=MatchLowercase}
  \defaultfontfeatures[\rmfamily]{Ligatures=TeX,Scale=1}
\fi
\IfFileExists{upquote.sty}{\usepackage{upquote}}{}
\IfFileExists{microtype.sty}{
  \usepackage[]{microtype}
  \UseMicrotypeSet[protrusion]{basicmath} 
}{}
\makeatletter
\@ifundefined{KOMAClassName}{
  \IfFileExists{parskip.sty}{%
    \usepackage{parskip}
  }{
    \setlength{\parindent}{0pt}
    \setlength{\parskip}{6pt plus 2pt minus 1pt}}
}{
  \KOMAoptions{parskip=half}}
\makeatother
\usepackage{xcolor}
\IfFileExists{xurl.sty}{\usepackage{xurl}}{} 
\IfFileExists{bookmark.sty}{\usepackage{bookmark}}{\usepackage{hyperref}}
\hypersetup{
  pdftitle={Abstract Orientable Incidence Structure and Algorithms for Finite Bounded Acyclic Categories. I. Incidence Structure},
  pdfauthor={Yu-Wei Huang},
  hidelinks,
  pdfcreator={LaTeX via pandoc}}
\usepackage{graphicx}
\makeatletter
\def\maxwidth{\ifdim\Gin@nat@width>\linewidth\linewidth\else\Gin@nat@width\fi}
\def\maxheight{\ifdim\Gin@nat@height>\textheight\textheight\else\Gin@nat@height\fi}
\makeatother
\setkeys{Gin}{width=\maxwidth,height=\maxheight,keepaspectratio}
\makeatletter
\def\fps@figure{htbp}
\makeatother
\providecommand{\tightlist}{%
  \setlength{\itemsep}{0pt}\setlength{\parskip}{0pt}}
\usepackage{fullpage}
\usepackage[]{biblatex}
\title{Abstract Orientable Incidence Structure and Algorithms for Finite
Bounded Acyclic Categories. I. Incidence Structure}
\author{Yu-Wei Huang\thanks{l28071504@gs.ncku.edu.tw}}
\date{\today}

\begin{document}
\maketitle
\begin{abstract}
A generalization of incidence relations in abstract polytope has been
explored, and parameterized surfaces are used as primers. The abstract
orientable incidence structure is defined as an algebraic model of
incidence relations, in which some algebraic properties in abtract
polytope theory are generalized. The geometric interpretation of
abstract orientable incidence structure are also discussed. The
orientable incidence structure in a semi-regular normal CW complex are
briefly investigated.
\end{abstract}

\hypertarget{introduction}{%
\section{Introduction}\label{introduction}}

The concept of incidence relations can be traced back to the face
lattice of classical polytopes, which is a combinatorial structure among
facets \autocite{mcmullen2002}. This structure can be studied as an
abstract polytope, which is a binary relation describing whether a facet
is contained by another facet. However, abstract polytope is limited in
its ability to describe certain types of incidence relations. For
instance, the convexity of polytopes ensures that the relation between
two facets is unique, leading to a poset structure in the abstract
polytope. The flatness of geometric objects also makes the shape
uniquely determined by the boundaries. By relaxing these constraints,
more fundamental structures of incidence relations can be revealed. Some
geometric complexes like simplicial complex and CW complex can be seen
as a generalization of polytopes \autocite{kozlov2008,bjorner1984}.
However, these approaches require the introduction of auxiliary
geometric objects and are not well-suited for studying incidence
structures. Although it is well known that incidence relations can be
manifested by posets, there has been little discussion of generalizing
this concept to acyclic categories in the same manner. Therefore, this
article will re-describe incidence relations between geometric objects,
and show that it has the structure of acyclic category.

\hypertarget{incidence-structure}{%
\section{Incidence Structure}\label{incidence-structure}}

\hypertarget{facet-and-incidence-relation}{%
\subsection{Facet and Incidence
Relation}\label{facet-and-incidence-relation}}

\begin{figure}
\hypertarget{fig:spherical-shell}{%
\centering
\includegraphics[width=0.4\textwidth,height=\textheight]{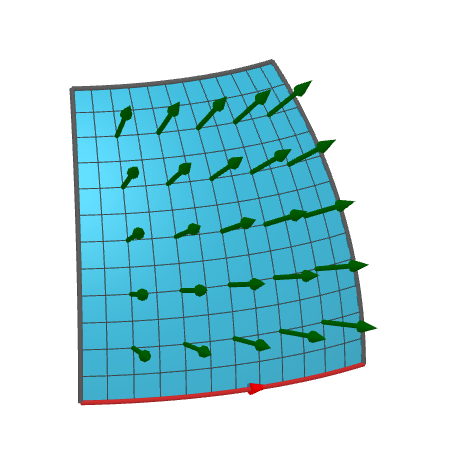}
\caption{A fregment of spherical shell. The arrows on the surface
represent its orientation. The arc on its bottom boundary, whose
orientation is represented as the arrow on the line, is
positively-oriented with respect to the spherical
shell.}\label{fig:spherical-shell}
}
\end{figure}

Below we will introduce incidence structure via parameterized
\(n\)-dimensional surfaces, which are referred to as \textbf{facets}.
Roughly speaking, a parameterized \(n\)-dimensional surface is a surface
defined by a function of \(n\) parameters. In order to rule out some
strange geometries, some constraints are adopted. The points on the
surface are described by the parameters called \textbf{coordinates}
faithfully, so that the topological properties of this surface are fully
described by the coordinate space. The coordinate space is not limit to
an Euclidean space, so that it can describe non-trivial topology like
torus. The values of coordinates are restricted to a given range in the
coordinate space, which should be a compact and connected open subset.
For example, the Figure \ref{fig:spherical-shell} is a facet of a
fragment of spherical shell described by
\(\{ S(\theta, \phi) \mid 0 < \theta < \pi/4, 0 < \phi < \pi/4 \}\),
where
\(S(\theta, \phi) \equiv \cos\theta \cos\phi \, \hat{x} + \sin\theta \cos\phi \, \hat{y} + \sin\phi \, \hat{z}\)
is a point on the surface with coordinate \((\theta, \phi)\). The
topology of this surface can be described by the coordinate range
\((0, \pi/4) \times (0, \pi/4) \subset \mathbb{R}^2\). To describe a
full sphere, the coordinate space cannot be an Euclidean space because
it cannot faithfully refer to all points on the sphere. It is useful for
doing actual geometric calculations if the coordinate space has a
metric, which is important when it is not a Cartesian coordinate system.
A good coordinate space should have small distortion, such that
geometric calculations can be performed accurately. A cone can be
described by \(\{ Y(x, y) | 0 < x^2 + y^2 < 1 \}\) with
\(Y(x, y) \equiv x \, \hat{x} + y \, \hat{y} + \, \sqrt{x^2 + y^2} \hat{z}\),
where the origin of the coordinate space \((0, 0)\) has been removed
since such point is not differentiable. This coordinate range is valid
although it is not homeomorphic to a disk; there is no need to cut it
off to make it homeomorphic to the disk. These examples shows how the
flexibility of coordinate helps us to deal with real geometry more
easily. As a special case, a point, described by a singleton set
\(\{ P \}\), has a coordinate space with zero parameter, which is just a
coordinate space that has only one valid coordinate.

If the coordinate space is orientable, the orientation of a facet can be
defined. For convenience, we assume they are always orientable. The
orientation of a facet at the point \(F(u, v, \dots)\), denoted as
\(\hat{d}F(u, v, \dots)\), is defined as the normalization of the wedge
product of the partial derivatives
\(dF(u, v, \dots) = \frac{\partial}{\partial u} F(u, v, \dots) \wedge \frac{\partial}{\partial v} F(u, v, \dots) \wedge \dots\).
In the example of Figure \ref{fig:spherical-shell}, the orientation is
\(\hat{d}S(\theta, \phi) = \cos\theta \cos\phi \, \hat{y} \wedge \hat{z} + \sin\theta \cos\phi \, \hat{z} \wedge \hat{x} + \sin\phi \, \hat{x} \wedge \hat{y}\),
which can be seen as a radial pseudovector. Each point in a
\(n\)-dimensional facet has an orientation as a normalized multivector
of degree \(n\) (more precisely, a \(n\)-blade). In the same sense, the
orientation of a line is just a multivector of degree 1, which is an
unit tangent vector; the orientation of a point is just a multivector of
degree 0, which only can be a sign. The orientation of a point cannot be
derived naturally, so one should attach a sign to define its
orientation, which is called a \textbf{signed point}.

A valid facet should be possible to be extended by including the
boundary of the coordinate range, which is called \textbf{closed
facets}. The image of the open coordinate range is called the
\textbf{inner} of the facet, and the image of the boundaries of the
coordinate range is called the \textbf{boundaries} of the facet. The
inner and the boundaries of a facet should be disjoint. Unlike the inner
of the facet, the boundaries are not faithfully described by
coordinates, but are at least locally faithfully. To describe the
boundaries of a given facet, a set of disjoint facets are introduced
such that the union of them is equal to the boundaries, and these facets
are called \textbf{subfacets} of the given facet. The relations between
a facet and its subfacets are called \textbf{incidence relations}, which
are described by embed function and local connectedness. \textbf{Embed
function} is an mapping between coordinate spaces of facets, which shows
how a facet is embedded into another one. For example, an arc
\(\{ C(\theta) \mid 0 < \theta < \pi/4 \}\) with
\(C(\theta) \equiv \cos\theta \, \hat{x} + \sin\theta \, \hat{y}\) is
covered by the surface of the example of Figure
\ref{fig:spherical-shell}, which can be described by the equation
\(C(\theta) = S(\theta, 0)\). This embed function is denoted as
\(\phi : C \to S\), where \(C\) is a subfacet of the facet \(S\). Notice
that there can be multiple embed functions between two facets. For
example, the circle, described by
\(\{ C(\theta) | 0 < \theta < 2\pi \}\) with
\(C(\theta) \equiv \cos\theta \, \hat{x} + \sin\theta \, \hat{y}\),
covers the point \(P \equiv 1 \, \hat{x} + 0 \, \hat{y}\) in two ways:
\(P = C(0)\) and \(P = C(2\pi)\). This is different from the
conventional polytope, in which there is only one relation between
incident facets.

The definiton of subfacets also implies two special facets. The null
face \(\varnothing\) is defined as an empty set, which has no valid
coordinate. The null face is usually considered as a
\((-1)\)-dimensional facet. The null face \(\varnothing\) is covered by
all facets, because there is a map from an empty set to any set.
Consider multiple disjoint facets with shared subfacets, the space
containing these facets, called the universe \(\mathbb{U}\), can also be
considered as a facet. For example, 3D Euclidean space is defined as
\(\mathbb{U}(x, y, z) = x \, \hat{x} + y \, \hat{y} + z \, \hat{z}\).
The universe \(\mathbb{U}\) covers all facets because it is the codomain
of functions defining facets. The null face \(\varnothing\) and the
universe \(\mathbb{U}\) are called improper facets.

\begin{figure}
\hypertarget{fig:orientation}{%
\centering
\includegraphics[width=0.7\textwidth,height=\textheight]{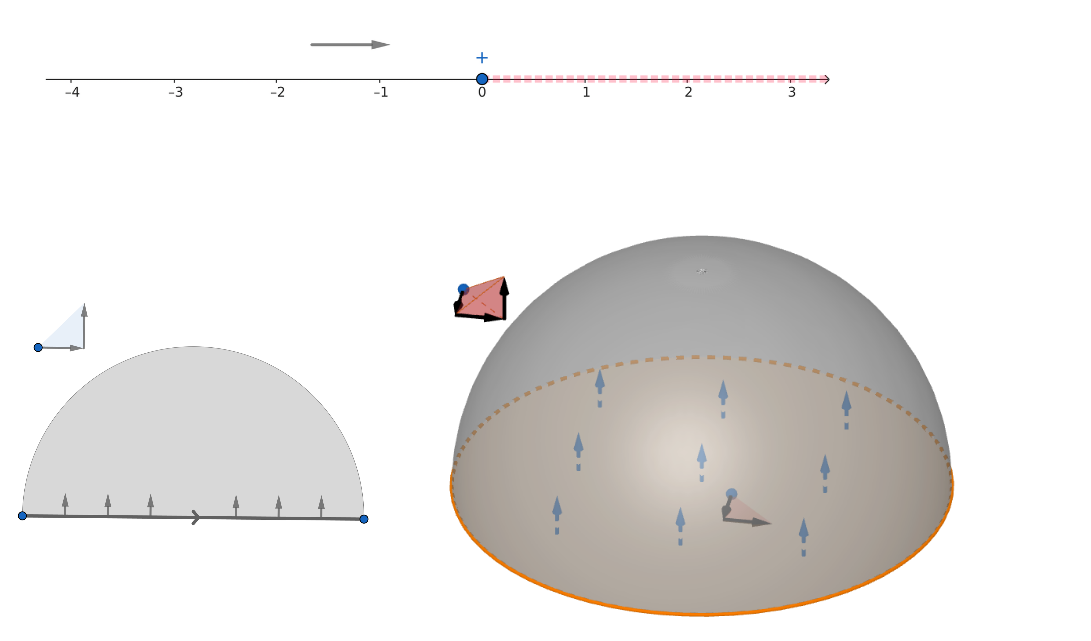}
\caption{The relative orientation of signed points, lines, and planes.
Their orientations are represented as a chain of arrows (drawn as
simplexes), so that the wedge product of these vectors forms the
orientation. The shaded regions represent positive side of
facets.}\label{fig:orientation}
}
\end{figure}

Another property of an incidence relation called the \textbf{local
connectedness} is manifested by the relative position between two
incident facets. If their dimensions differ by only \(1\), the relative
position can always be described by a sign (denoted as \(|\phi|\)),
determining which the side of the subfacet faces the inner of the facet.
In the example of Figure \ref{fig:spherical-shell}, the spherical shell
is on the positive side of the arc according to the right-hand rule, so
we define \(|\phi| = +1\) (see Figure \ref{fig:spherical-shell}).
Generalizing to incidence relation between \(n\)-dimensional and
\((n{+}1)\)-dimensional facets, we define the orientation between them
as the difference of wedge products: if the \(n\)-dimentional subfacet
\(B\) is embedding into the \((n{+}1)\)-dimentional facet \(F\) by
mapping \(\phi : B \to F\), define \(|\phi|\) as a sign such that
\(\hat{d}B \wedge \hat{n} = |\phi| \hat{d}F\), where \(\hat{d}F\) and
\(\hat{d}B\) are their orientations at a point in the subfacet \(B\),
and \(\hat{n}\) is a vector on this point pointing to the inner of the
facet \(F\). \(|\phi| = +1\) (\(|\phi| = -1\)) means the subfacet \(B\)
is positively- (negatively-) oriented with respect to the facet \(F\)
(see Figure \ref{fig:orientation}). Note that all points in the subfacet
\(B\) should give the same value of \(|\phi|\). Under this definition,
one can say: a line segment has one positive-oriented subfacet and one
negative-oriented subfacet, and both of them are positive points at the
endpoints; a circle with counter-clockwise tangent vectors is a
positive-oriented subfacet of the enclosed disk; a sphere defined like
the above example is negatively-oriented with respect to the enclosed
ball. The incidence relation between the null face and a signed point is
a special case: positive (negative) point always is positively-
(negatively-) oriented with respect to the null face. A subfacet can be
covered by a facet with multiple orientations. For example, consider a
plane \(S\) and a line segment \(L\) which is in the middle of the
plane. Since the plane is on the left and right of this line, there are
two choices of the inward pointing vectors. Then the incidence relations
between them are described by \(\phi_+\) and \(\phi_-\), which are the
two identical embed functions with different orientations by letting
\(|\phi_+| = +1\) and \(|\phi_-| = -1\).

\begin{figure}
\hypertarget{fig:diamond}{%
\centering
\includegraphics[width=0.5\textwidth,height=\textheight]{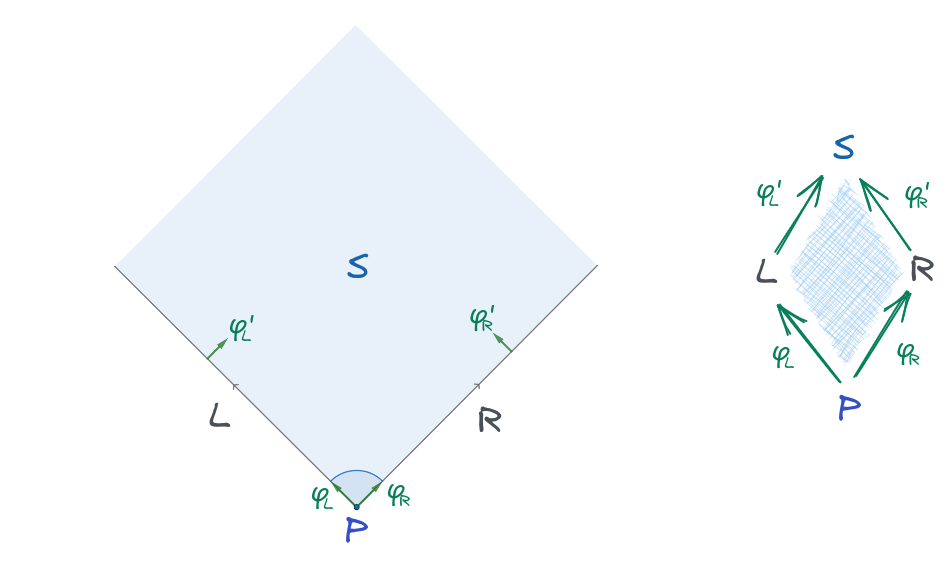}
\caption{The incidence relations between a square and its vertex. The
incidence relations from vertex to edge and from edge to square can be
composed, resulting in the incidence relation from the vertex to the
square. These incidence relations form a diamond shape (right diagram),
which is filled up to indicate the equivalence of composition. The
equivalence of composition represents the completeness of the vertex,
drawn as a fan shape on the vertex.}\label{fig:diamond}
}
\end{figure}

\begin{figure}
\hypertarget{fig:crescent}{%
\centering
\includegraphics[width=0.5\textwidth,height=\textheight]{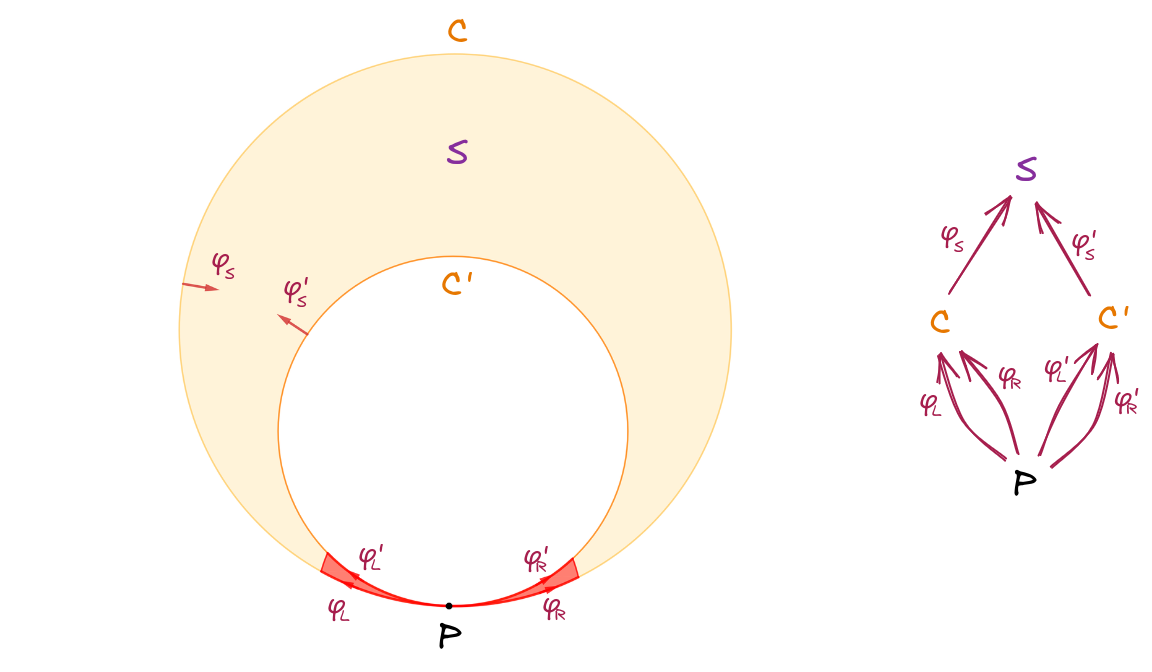}
\caption{Multiple incidence relations on a crescent shape. There are two
diamonds on the point \(P\), which represents two directions (drawn as
two red fan shapes) go from this point to the inner of the crescent
shape. The right diagram represents the embed functions between facets,
while it's difficult to be filled up to represent their diamond
properties.}\label{fig:crescent}
}
\end{figure}

It is obviously that incidence relations are composable: if a point can
embed into a line, and the line can embed into another surface, then
this point can also embed into this surface. The equivalence relations
between composed embed functions manifest the local connectedness, which
is not just a sign. The incidence relations and their composition rules
constitute an \textbf{incidence structure}. Figure \ref{fig:diamond}
shows the incidence relations between a square and its vertex. The
vertex point embeds into a square by mapping \(\phi\), which is equal to
composition \(\phi_L' \circ \phi_L\) and \(\phi_R' \circ \phi_R\), where
\(\phi_L\) (\(\phi_R\)) is the incidence relation from this point to the
left (right) adjacent edge, and \(\phi_L'\) (\(\phi_R'\)) is the
incidence relation from the left (right) adjacent edge to the square.
Such equivalence relation is just the diamond property in abstract
polytope, and this property describes the completeness of the vertex. In
abstract polytope, this diamond always exists and is unique for any two
incident facets with dimension differing by \(2\). If there is a shape
that has two vertices sharing the same point, this point is covered by
this area in two directions. For example, consider the area between two
internally tangent circles called crescent shape (see Figure
\ref{fig:crescent}), the tangent point is covered by outer circle in two
directions, say \(\phi_L\) and \(\phi_R\), and outer circle is also
covered by the crescent shape, say \(\phi_S\). Then one can go from the
tangent point to the inner of the crescent shape in two ways, which is
equal to \(\phi_S \circ \phi_L\) and \(\phi_S \circ \phi_R\). These two
incidence relations are not equivalent, so they don't form a diamond;
the angle between these two edges has been separated by the inner
circle. Instead, the diamonds are formed by
\(\phi_S \circ \phi_L = \phi_S' \circ \phi_L'\) and
\(\phi_S \circ \phi_R = \phi_S' \circ \phi_R'\), where \(\phi_L'\),
\(\phi_R'\) and \(\phi_S'\) are the same as above but relative to the
inner circle. These two diamonds represent two different directions of
incidence relations between the tangent point and the crescent shape;
there are two separated ranges of angles that can go from this point to
the inner. It shows the difference between diamonds indicates the local
connectedness around the point.

Not all incidence relations can be defined continuously. The Riemann
surface \(R(z) \equiv (z, f(z))\) for the function
\(f(z) = \sqrt{z} (1 - |z|)\) with \(0 < |z| < 1\) is a
\(2\)-dimensional facet in a \(4\)-dimensional space, which doubly cover
the unit circle \(C(z) \equiv (z, 0)\) with \(|z| = 1\). Their are two
incidence relations: any point on the circle, say \(C(e^{i\theta})\),
should be mapped to \(R(e^{i\theta})\) and \(R(e^{i(2\pi + \theta)})\),
representing two embed functions respectively. Depending on the choice
of the branch cut, the embed functions are discontinuous at that point.
It is impossible to write such embed functions as continuous functions.
Breaking of continuity may cause some geometric calculation problems. To
solve this problem, the calculation should be limited to a small range,
such that in this range the embed functions can be written as continuous
functions depending on where it is. A facet with a non-orientable
subfacet will also encounter similar problems, in which the incidence
relations cannot have a continuous sign. These kind of incidence
relations are said to be \textbf{non-orientable}. To simplify the
discussion, we only focus on the incidence structure with orientable
incidence relations, which we call an \textbf{orientable incidence
structure}.

\hypertarget{vertex-figure-and-local-connectedness}{%
\subsection{Vertex Figure and Local
Connectedness}\label{vertex-figure-and-local-connectedness}}

\begin{figure}
\hypertarget{fig:vertex-figure}{%
\centering
\includegraphics[width=0.4\textwidth,height=\textheight]{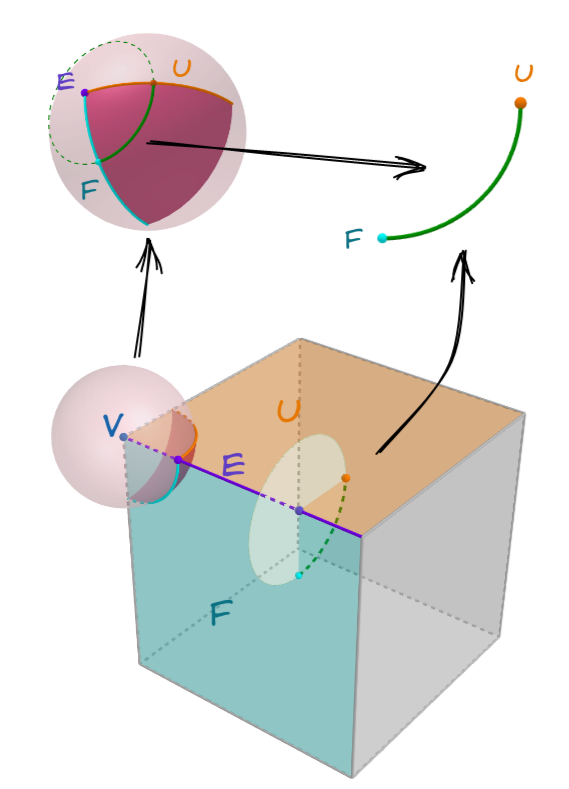}
\caption{A vertex figure (top left) and edge figure (top right) of a
cube (bottom). Where \(U\) and \(F\) denote two faces, \(E\) denotes an
edge, and \(V\) denote a vertex. The vertex figure on the vertex \(V\)
is formed by sliced images of those facets. Further taking vertex figure
on the vertex \(E\) is the same as taking vertex figure on the edge
\(E\) in the original cube.}\label{fig:vertex-figure}
}
\end{figure}

To explain the concept of local connectedness more clearly, we shall
introduce vertex figures, which is defined in the same manner as
polytope. To make a vertex figure on a \(0\)-dimensional facet (a
positive point), one puts a ball \(D^n\) on it, which is small enough so
that only the facets that cover this point intersect this ball. Then the
sliced images of facets (the intersection of the boundary of the ball
\(\partial D^n\) and facets) are defined as the \textbf{vertex figure}
of the corresponding facets. For example, the vertex figure on a vertex
of a cube is a spherical triangle, which is an intersection of the cube
and a sphere centered at this vertex (see Figure
\ref{fig:vertex-figure}). The orientation of the sliced image \(L\) of
the facet \(F\) should obey \(\hat{n} \wedge \hat{d}L = \hat{d}F\),
where \(\hat{n}\) is outward pointing vector from the center of the
slicing ball. The incidence relations of this vertex figure can be
trivially derived by restricting the domain and codomain of mapping.
It's easy to prove that the sign of incidence relations are the same.
More general, the vertex figure on any facet can be defined as: pick a
point on given \(m\)-dimensional facet \(V\), and put a ball \(D^{n-m}\)
on this point, which is perpendicular to its orientation \(\hat{d}V\) at
this point. Then the sliced images of facets (the intersections of the
boundary of this ball \(\partial D^{n-m}\) and facets) are defined as
the vertex figure of the corresponding facets on the facet \(V\). All
points on the facet \(V\) should give the same vertex figure. The
orientation of the sliced image \(L\) of the facet \(F\) should obey
\(\hat{d}V \wedge \hat{n} \wedge \hat{d}L = \hat{d}F\), where
\(\hat{n}\) is outward pointing vector from the center of the the
slicing ball. For example, the vertex figure on an edge (i.e., the edge
figure) of a cube is a segment of arc, which is the intersection of the
cube and a circle with this edge as its axis (see Figure
\ref{fig:vertex-figure}). Its obviously the operations of taking vertex
figure are composable. For example, the edge figure of a cube is the
vertex figure of the vertex figure, as shown in Figure
\ref{fig:vertex-figure}.

Some sliced images on a subfacet have multiple connected components,
each of them is defined as an individual facet in this vertex figure.
This kind of facets are said to be \textbf{locally separated} by this
subfacet. For example, the vertex figure of a plane on a point lying in
the middle of the plane is a circle, which is connected, so we say the
plane is not locally separated by this point. But the vertex figure of a
plane on a line lying in the middle of the plane becomes two points,
which have two connected components, so we say the plane is locally
separated by this line. The crescent shape (Figure \ref{fig:crescent})
is a non-trivial example, whose vertex figures on the tangent point are
two line segments, so we say the crescent shape is locally separated by
the tangent point. It shows the local connectedness property described
by incidence relation becomes connectedness property of the vertex
figure, and such correspondence is naturally required if one considers
the incidence structure of the vertex figure. Also, the local
connectedness property can be propagated via connectivity, so incidence
relations can be tested in every points in the connected subfacet. These
descriptions manifest the local and global properties of incidence
relations.

\begin{figure}
\hypertarget{fig:torus}{%
\centering
\includegraphics[width=0.6\textwidth,height=\textheight]{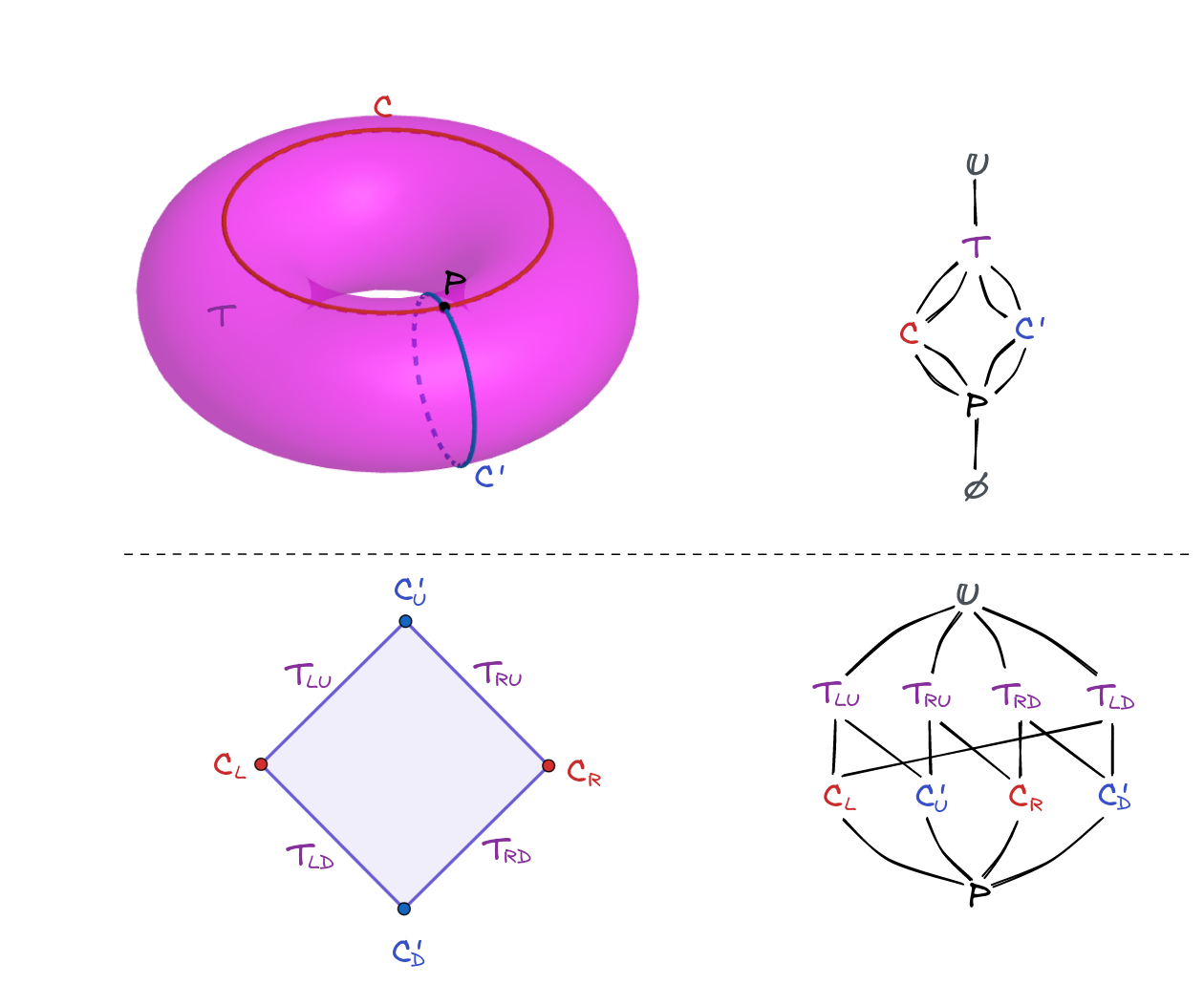}
\caption{The vertex figure of a torus. The top is the torus (left) and
its Hasse diagram (right); the bottom is the vertex figure on point
\(P\) (left) and the corresponding Hasse diagram (right), which is much
more complicated.}\label{fig:torus}
}
\end{figure}

The orientable incidence structure can be drawn as directed acyclic
multigraphs, where nodes are facets, and edges are nondecomposable embed
functions. In this diagram, called the \textbf{Hasse diagram}, all
incident facets can be connected by directed paths, manifesting the
poset structure of incidence relations. The Hasse diagram is not enough
to describe an orientable incidence structure, since the composition
rules for incidence relations are not shown. Because sliced images of
connected facets may become disconnected, the incidence structure of a
vertex figure is not just a subgraph of the Hasse diagram. For example,
the vertex figure on the intersection point of two circles on a torus is
much more complicated to itself (see Figure \ref{fig:torus}), but the
composition rules of incidence relations are always the same.

In our theory, the incidence relations are not unique to two incident
facets. The difference of two incidence relations represents the local
connectedness of the facet, which is the main difference to the
conventional incidence structure. The existence of multiple incidence
relations forces us to consider the composition rules between them.
Unlike abstract polytope, where the structure can be described by a
poset, the orientable incidence structure act like a category. We will
extract out such structure in the next section.

\hypertarget{abstract-orientable-incidence-structure}{%
\section{Abstract Orientable Incidence
Structure}\label{abstract-orientable-incidence-structure}}

\hypertarget{bounded-acyclic-category-and-bounded-poset}{%
\subsection{Bounded Acyclic Category and Bounded
Poset}\label{bounded-acyclic-category-and-bounded-poset}}

The abstract orientable incidence structure is an algebraic structure
which represents the orientable incidence relations between facets.
Unlike previous section, it only captures its combinatorial nature
without specifying geometric properties. An \textbf{abstract orientable
incidence structure}, denoted as \(\mathcal{F}\), is a \textbf{graded
bounded acyclic category} satisfying \textbf{semi-diamond property}. An
\textbf{acyclic category} is a category without non-trivial cycles. It
is impossible to have two non-identity morphisms in opposite directions
at the same time, so morphisms can only go in one direction. For a
finite acyclic category, there is a rank function that maps an object to
an integer, so that all morphisms increase the rank of object except
identity. One can define the rank of an object as the dimension of
corresponding facet. The acyclic category with a rank function is called
\textbf{graded acyclic category}. \textbf{Bounded} means it has exactly
one initial object and terminal object. The initial object, which has
the lowest rank, represents the null face, denoted as \(\varnothing\);
the terminal object, which has the highest rank, represents the
universe, denoted as \(\mathbb{U}\). Morphisms represent incidence
relations between the facets. The morphism \(\phi : F \rightarrow S\)
means that facet \(F\) is covered by facet \(S\). There may have
multiple morphisms between two objects, and this property represents the
local connectedness of the incidence relation, as we discussed above.
Below we mainly focus on the property of bounded acyclic category, and
compare it with the bounded poset.

Recall the definition of a \textbf{bounded poset} (a poset which
includes its lower and upper bounds), which is the base of abstract
polytope:

\begin{enumerate}
\def\labelenumi{\arabic{enumi}.}
\tightlist
\item
  \emph{irreflexivity}: \(a \nless a\).
\item
  \emph{asymmetry}: if \(a < b\) then \(b \nless a\).
\item
  \emph{transitivity}: if \(a < b\) and \(b < c\) then \(a < c\).
\item
  \emph{least}: there exists exactly one \(\varnothing\) such that
  forall \(a \neq \varnothing\) there is \(\varnothing < a\).
\item
  \emph{greatest}: there exists exactly one \(\mathbb{U}\) such that
  forall \(a \neq \mathbb{U}\) there is \(a < \mathbb{U}\).
\end{enumerate}

Now compare to the definition of \textbf{bounded acyclic category}:

\begin{enumerate}
\def\labelenumi{\arabic{enumi}.}
\tightlist
\item
  \emph{irreflexivity}: there is no \(\phi : F \rightarrow F\) for any
  \(\phi \neq \operatorname{id}\).
\item
  \emph{asymmetry}: if there is \(\phi : F \rightarrow S\) for some
  \(\phi \neq \operatorname{id}\), then there is no
  \(\phi' : S \rightarrow F\) for any \(\phi' \neq \operatorname{id}\).
\item
  \emph{transitivity}: if there are \(\phi : P \rightarrow F\) and
  \(\phi' : F \rightarrow S\) for some \(\phi \neq \operatorname{id}\)
  and \(\phi' \neq \operatorname{id}\), then
  \(\phi' \circ \phi : P \rightarrow S\).
\item
  \emph{initial}: there exists exactly one \(\varnothing\), such that
  forall \(F \neq \varnothing\) there is
  \(\phi : \varnothing \rightarrow F\) for \emph{exactly one}
  \(\phi \neq \operatorname{id}\).
\item
  \emph{terminal}: there exists exactly one \(\mathbb{U}\), such that
  forall \(F \neq \mathbb{U}\) there is
  \(\phi : F \rightarrow \mathbb{U}\) for \emph{exactly one}
  \(\phi \neq \operatorname{id}\).
\end{enumerate}

The key different is that, there can be multiple relations between two
objects in the bounded acyclic category, and the transitivity of order
relations becomes the composition of morphisms. A bounded poset is just
a thin bounded acyclic category. It is natural to define an induced
bounded poset for a bounded acyclic category by forgetting differences
between multiple relations. In this sense, one says acyclic category is
a generalization of poset. Abstract orientable incidence structure also
loosens three constraints of abstract polytope, which will be discussed
in the end of this section. Notice that the rules \emph{least} and
\emph{greatest} in a bounded poset becomes \emph{initial} and
\emph{terminal}, which looks stricter in the bounded acyclic category.
The initial object and the terminal object are called improper objects.

\hypertarget{upper-category-and-downward-functor}{%
\subsection{Upper Category and Downward
Functor}\label{upper-category-and-downward-functor}}

The upper closure of an element \(x\) in a bounded poset \(\mathcal{P}\)
is a subset whose element is greater than or equal to \(x\). Such subset
is also a bounded poset, where the order relation is inherited from the
host poset \(\mathcal{P}\). To emphasize this, it is denoted as
\(\mathcal{P} \operatorname{\uparrow} x\). Similarly, one can construct
a full subcategory \(\mathcal{F}'\) by only including the upper closure
of an object \(F_m\). There is no non-trivial morphism that point to the
object \(F_m\), but it is not an initial object since the morphism
\(\phi : F_m \to S\) may not be unique for the object \(S\), so that it
is not a bounded acyclic category. To fix it, one should split objects
such that there is only one initial morphism for each object. Construct
\textbf{upper category} of object \(F_m\) in a bounded acyclic category
\(\mathcal{F}\), denoted as \(\mathcal{F} \operatorname{\uparrow} F_m\):
\(\operatorname{Obj}(\mathcal{F} \operatorname{\uparrow} F_m)\) is a set
of \(F_s | \phi_{sm} \rangle\) for all objects \(F_s\) and morphisms
\(\phi_{sm} : F_m \rightarrow F_s\).
\(\operatorname{Hom}(F_s | \phi_{sm} \rangle, F_t | \phi_{tm} \rangle)\)
is a set of
\(\phi_{ts} | \phi_{sm} \rangle : F_s | \phi_{sm} \rangle \rightarrow F_t | \phi_{tm} \rangle\)
for all morphisms \(\phi_{ts} : F_s \rightarrow F_t\) and
\(\phi_{sm} : F_m \to F_s\) such that
\(\phi_{ts} \circ \phi_{sm} = \phi_{tm}\). The composition rule is
\(\phi_{pt} | \phi_{tm} \rangle \circ \phi_{ts} | \phi_{sm} \rangle = \phi_{ps} | \phi_{sm} \rangle\)
with \(\phi_{ps} = \phi_{pt} \circ \phi_{ts}\). The object \(F_t\) have
been marked with morphism \(\phi_{tm}\), so that it is splitted into
multiple objects corresponding to each initial morphism \(\phi_{tm}\).
Although the object splits, the composition rule is inherited from the
host category \(\mathcal{F}\). In terms of incidence structure, the
upper category \(\mathcal{F} \operatorname{\uparrow} F_m\) just
corresponds to the vertex figure on an object \(F_m\): in upper
category, an object \(F_n | \phi_{nm} \rangle\) represents a connected
component of sliced images, and a morphism
\(\phi_{ts} | \phi_{sm} \rangle\) represents an incidence relation
between facets \(F_s | \phi_{sm} \rangle\) and
\(F_t | \phi_{tm} \rangle\).

There is a functor
\(F_m^{\downarrow} : \mathcal{F} \operatorname{\uparrow} F_m \rightarrow \mathcal{F}\)
that maps object \(F_s | \phi_{sm} \rangle\) to \(F_s\), and maps
morphism
\(\phi_{ts} | \phi_{sm} \rangle : F_s | \phi_{sm} \rangle \rightarrow F_t | \phi_{tm} \rangle\)
to \(\phi_{ts} : F_s \rightarrow F_t\). This functor is a reversed
operation of taking upper category, so it is named \textbf{downward
functor}. Drawn as Hasse diagrams, the downward functor becomes graph
homomorphism between two graphs. Unlike the abstract polytope, this is
not injective on nodes and edges, which manifests the local
connectedness of incidence relations.

\begin{figure}
\hypertarget{fig:upper-categories}{%
\centering
\includegraphics[width=1\textwidth,height=\textheight]{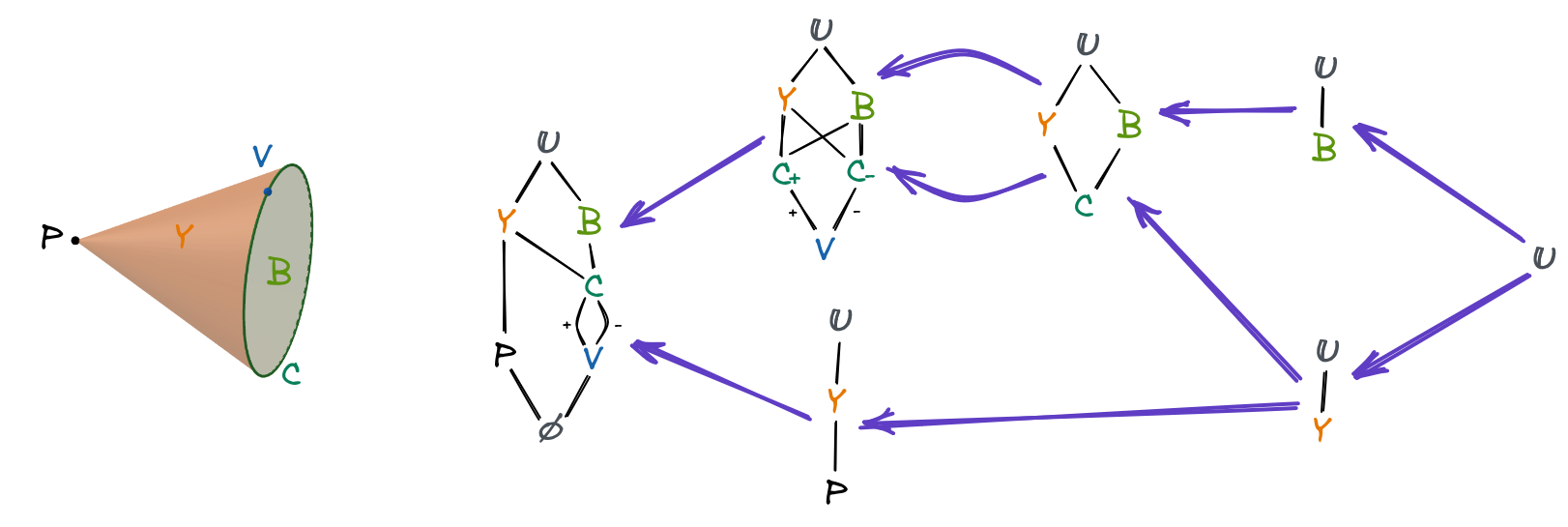}
\caption{The category of upper categories (right) of a cone (left). It
is drawn as a set of Hasse diagrams with graph homomorphisms between
them (purple arrows). This representation is sufficiently to describe a
bounded acyclic category.}\label{fig:upper-categories}
}
\end{figure}

Since an upper category itself is also a bounded acyclic category, one
can further construct an upper category of an object in the upper
category. The category
\(\mathcal{F} \operatorname{\uparrow} F_m \operatorname{\uparrow} F_n | \phi_{nm} \rangle\)
is an upper category of object \(F_n | \phi_{nm} \rangle\), where
\(\operatorname{Obj}(\mathcal{F} \operatorname{\uparrow} F_m \operatorname{\uparrow} F_n | \phi_{nm} \rangle)\)
is a set of object
\(F_s | \phi_{sm} \rangle | \phi_{sn} | \phi_{nm} \rangle \rangle\), or
abbreviated as \(F_s | \phi_{sn}, \phi_{nm} \rangle\), for all objects
\(F_s\) and morphisms \(\phi_{sn} : F_n \rightarrow F_s\);
\(\operatorname{Hom}(F_t | \phi_{tn}, \phi_{nm} \rangle, F_s | \phi_{sn}, \phi_{nm} \rangle)\)
is a set of
\(\phi_{ts} | \phi_{sm} \rangle | \phi_{sn} | \phi_{nm} \rangle \rangle : F_s | \phi_{sm} \rangle | \phi_{sn} | \phi_{nm} \rangle \rangle \rightarrow F_t | \phi_{tm} \rangle | \phi_{tn} | \phi_{nm} \rangle \rangle\),
or abbreviated as
\(\phi_{ts} | \phi_{sn}, \phi_{nm} \rangle : F_s | \phi_{sn}, \phi_{nm} \rangle \rightarrow F_t | \phi_{tn}, \phi_{nm} \rangle\)
for all morphisms \(\phi_{ts} : F_s \rightarrow F_t\) and
\(\phi_{sn} : F_n \rightarrow F_s\) such that
\(\phi_{ts} \circ \phi_{sn} = \phi_{tn}\). Noting that we assume
\(\phi_{sm} = \phi_{sn} \circ \phi_{nm}\) and
\(\phi_{tm} = \phi_{tn} \circ \phi_{nm}\), and the abbreviation works
because \(\phi_{sm}\) can be determined uniquely. It is isomorphic to
the upper category \(\mathcal{F} \operatorname{\uparrow} F_n\), which
shows that taking upper category is composable. This rule is consistent
with above discussion: taking vertex figure is composable. Similarly,
downward functors are also composable. Downward functor can be reduced
down to a functor
\(\phi_{mn}^{\downarrow} : \mathcal{F} \operatorname{\uparrow} F_m \rightarrow \mathcal{F} \operatorname{\uparrow} F_n\)
for morphism \(\phi_{mn} : F_n \rightarrow F_m\), in which object
\(F_s | \phi_{sm} \rangle\) is mapped to \(F_s | \phi_{sn} \rangle\),
and morphism
\(\phi_{ts} | \phi_{sm} \rangle : F_s | \phi_{sm} \rangle \rightarrow F_t | \phi_{tm} \rangle\)
is mapped to
\(\phi_{ts} | \phi_{sn} \rangle : F_s | \phi_{sn} \rangle \rightarrow F_t | \phi_{tn} \rangle\),
where \(\phi_{sn} = \phi_{sm} \circ \phi_{mn}\). Upper closures of a
poset form a poset by inclusion relation, its opposite category is
isomorphic to the host poset. Similarly, upper categories
\(\{ \mathcal{F} \operatorname{\uparrow} F_m | F_m \in \operatorname{Obj} \}\)
and downward functors
\(\{ \phi_{mn}^{\downarrow} | \phi_{mn} \in \operatorname{Hom} \}\) also
form a bounded acyclic category, its opposite category is isomorphic to
the host category \(\mathcal{F}\) (see Figure
\ref{fig:upper-categories}).

The Hasse diagram of a bounded acyclic category can easily see which
part is splittable. Two proper objects are said to be \textbf{linked} if
they are connected in the Hasse diagram excluding improper objects. If
not all proper objects are linked, this category is said to be
\textbf{splittable}. It is the generalization of connectedness condition
of abstract polytope. The linked objects may form multiple clusters, and
they can be splitted into separated bounded acyclic categories. The
linked clusters can be determined by linkages between
\textbf{minimal/maximal objects}, which are defined as minimal/maximal
elements of the induced poset of all proper objects. A bounded acyclic
category is said to be \textbf{strongly unsplittable} iff all sections
are unsplittable, which is an analogue of strongly connectedness of
posets. The definition of sections in an acyclic category will be
introduced latter. The words ``linked'' and ``splittable'' are used
because ``disjoint'' and ``disconnected'' are reserved for actual
geometric properties.

\hypertarget{section-category-and-morphism-chain}{%
\subsection{Section Category and Morphism
Chain}\label{section-category-and-morphism-chain}}

\begin{figure}
\hypertarget{fig:face-figure}{%
\centering
\includegraphics[width=0.5\textwidth,height=\textheight]{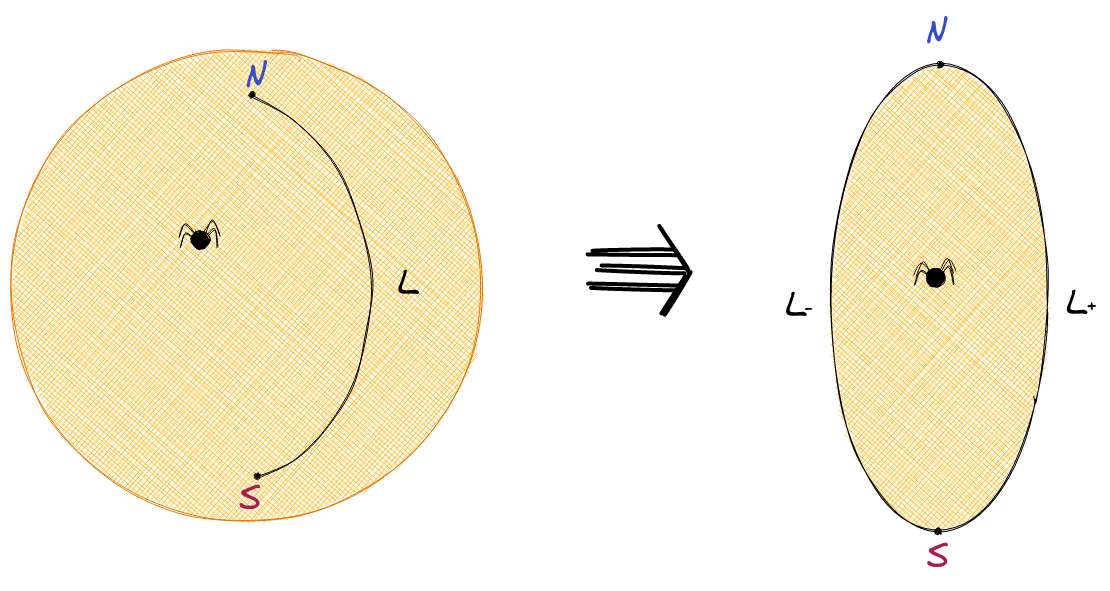}
\caption{A 2D ant lives on a sphere with a meridian (left), which is the
same as living on its face figure from the ant's point of view
(right).}\label{fig:face-figure}
}
\end{figure}

There is a dual concept of upper category called \textbf{lower
category}, which is a generalization of lower closure of a poset. Upper
category is constructed by treating one object as the initial object,
and split objects such that outgoing morphisms are initial morphisms. In
the opposite, lower category, denoted as
\(\mathcal{F} \operatorname{\downarrow} F_m\), is constructed by
treating one object as the terminal object, and split objects such that
incoming morphisms are terminal morphisms:
\(\operatorname{Obj}(\mathcal{F} \operatorname{\downarrow} F_m)\) is a
set of \(\langle \phi_{mt} | F_t\) for all objects \(F_t\) and morphisms
\(\phi_{mt} : F_t \rightarrow F_m\);
\(\operatorname{Hom}(\langle \phi_{ms} | F_s, \langle \phi_{mt} | F_t)\)
is a set of
\(\langle \phi_{mt} | \phi_{ts} : \langle \phi_{ms} | F_s \rightarrow \langle \phi_{mt} | F_t\)
for all morphisms \(\phi_{ts} : F_s \rightarrow F_t\) and
\(\phi_{mt} : F_t \rightarrow F_m\) such that
\(\phi_{mt} \circ \phi_{ts} = \phi_{ms}\). The composition rule is
\(\langle \phi_{mp} | \phi_{pt} \circ \langle \phi_{mt} | \phi_{ts} = \langle \phi_{mp} | \phi_{ps}\)
with \(\phi_{ps} = \phi_{pt} \circ \phi_{ts}\). In geometry, lower
category corresponds to \textbf{face figure}, which is like imagining an
ant living in a 2D world. The face figure of a sphere with one meridian
is equivalent to a 2D space bounded by two lines (see Figure
\ref{fig:face-figure}). The line is duplicated because it is covered by
the sphere in two ways. As a creature living in 3D space, we know it is
one line, but for the 2D creatures living on the sphere, they only see a
space rift, and they cannot confirm if there is a wider space in this
rift.

There is also a categorical analogue of the interval of a poset:
\([a, b] \equiv \{ c \in \mathcal{P} | a \le c \le b \}\), which can be
seen as a lower closure of an upper closure. Given a poset, by taking
non-trivial intervals between all minimal and maximal elements, the
given poset can be covered by multiple bounded posets. Similarly, an
acyclic category can be turned into direct sum of bounded acyclic
categories by splitting objects, and those bounded acyclic categories
are called \textbf{sections}, whose name comes from abstract polytope.
Consider an acyclic category, say \(\mathcal{F}\), define a
\textbf{section category} of morphism \(\phi_0\) on \(\mathcal{F}\),
denoted as \(\langle \mathcal{F} \rangle_{\phi_0}\):
\(\operatorname{Obj}(\langle \mathcal{F} \rangle_{\phi_0})\) is a set of
\(\langle \phi_n^\ast | \phi_n \rangle\) for all morphisms
\(\phi_n^\ast\) and \(\phi_n\) such that
\(\phi_n^\ast \circ \phi_n = \phi_0\).
\(\operatorname{Hom}(\langle \phi_n^\ast | \phi_n \rangle, \langle \phi_m^\ast | \phi_m \rangle)\)
is a set of \(\langle \phi_{m}^\ast | \phi_{mn} | \phi_{n} \rangle\) for
all morphisms \(\phi_{m}^\ast\), \(\phi_{mn}\), \(\phi_{n}\) such that
\(\phi_{m}^\ast \circ \phi_{mn} \circ \phi_{n} = \phi_0\).

taking section category doesn't break the algebra of morphisms, it just
renames some objects so that it is bounded. ``The algebra of morphisms''
is the generalization of ``the transitivity of order relation''. Define
\textbf{local-embedding function} for a poset: an order-preserving
function \(f\) is local-embedding iff it is an order isomorphism between
\([x, y]\) and \([f(x), f(y)]\) for any pair \(x \le y\). That is, the
element \(f(x) \le z' \le f(y)\) implies there exists exactly one
\(x \le z \le y\) such that \(f(z) = z'\). Note that it is not an
one-to-one function. The generalization to an acyclic category becomes a
\textbf{local-embedding functor}: for any morphism \(\phi\), the pair of
morphisms \(\psi'\) and \(\psi''\) with constraint
\(\psi'' \circ \psi' = \mu(\phi)\) implies there exists exactly one pair
of morphisms \(\phi'\) and \(\phi''\) with constraint
\(\phi'' \circ \phi' = \phi\) such that \(\mu(\phi') = \psi'\) and
\(\mu(\phi'') = \psi''\). The inverse of taking section categories is
the local-embedding functor. Local-embedding functor describes the
correspondence of algebra of morphisms without mention any object. It
shows that the identity of object isn't important in a bounded acyclic
category; initial morphisms and terminal morphisms are enough to
indicate objects in a bounded acyclic category.

The chain representation provides a clear view of this. Recall the
definition of the chain of a poset, which is defined as a total ordering
subset. A subset of a chain is also total ordering, called a subchain.
In an acyclic category, the morphism between two objects is not unique,
so the generalization becomes: the \textbf{morphism chain} of a morphism
\(\phi_{n0} : F_0 \rightarrow F_n\) is a non-empty sequence of morphisms
\(\langle \phi_n, \dots, \phi_2, \phi_1 \rangle\) that composed to this
morphism. Its \textbf{subchain} can be constructed by composing adjacent
morphisms into one. In other words, a subchain just skips some
intermediate objects except the start and end, which is slightly
different from the conventional definition of the morphism chain. A
morphism chain containing identity morphisms is said to be degenerated.
A morphism chain with \(n\) intermediate object is called \(n\)-chain.
\(0\)-chain is just the host morphism itself
\(\langle \phi_{n0} \rangle\), so it is said to be a trivial chain.
\(1\)-chain represents an object of the section of the morphism
\(\phi_{n0}\), and \(2\)-chain represents a morphism, whose source
object and target object are just two non-trivial subchains.

Abstract orientable incidence structure can be described by only their
subchain relations, called a \textbf{nerve}, which is also slightly
different from the conventional definition. The nerve of an acyclic
category is composed by a collection of \(n\)-chains
\(N_n = \{ \langle \phi_n, \dots, \phi_1, \phi_0 \rangle \}\) for
\(n = 0, 1, 2, \dots\) and the face maps \(d_i : N_n \to N_{n-1}\) for
\(i = 1 \sim n\), which skip the \(i\)-th object in the \(n\)-chain, and
the degeneracy maps \(s_i : N_n \to N_{n+1}\) for \(i = 1 \sim n\),
which insert an identity morphism at the \(i\)-th object. The laws for
face maps and degeneracy maps are the same as the conventional one,
which should not be repeated here. Note that in our definition, a
\(1\)-chain represents an object, not a morphism. A bounded acyclic
category has nerve with only one \(0\)-chain, which will back to the
conventional definition by replacing \(n\) with \(n-1\). The upper
category of a \(1\)-chain \(F\) becomes easy to define in this
formulation, which is just letting \(N_n' = d_1^{-n}(F)\) and
\(d_i' = d_i\), \(s_i' = s_i\), where \(d_1^{-n}\) indicates taking
preimage of \(d_1\) \(n\) times. Lower category and section category can
be defined in a similar approach. Two acyclic categories will have the
same nerve if they are surjectively section-embedded by the same direct
sum of bounded acyclic categories. Because no relabeling is required,
nerves are more natural for describing orientable incidence structures.

A section category not only can be defined on an acyclic category, but
also on the \textbf{non-recursive category}, which is a category without
non-trivial non-degenerated recursive morphism chain, that is,
\(\psi \circ \phi \circ \psi' = \phi\) implies \(\psi\) and \(\psi'\)
being identity morphisms. This property only rules out the infinity of
morphisms, and non-trivial cycles of morphisms are still possible, while
ensuring that their section categories are bounded acyclic categories.
Non-recursiveness and finiteness naturally induce acyclicity.

\hypertarget{semi-regular-normal-cw-complex}{%
\subsection{Semi-Regular Normal CW
Complex}\label{semi-regular-normal-cw-complex}}

The \textbf{diamond property} of an abstract orientable incidence
structure is stated as: a 2-rank morphism should be divided into
\emph{exactly} two maximum chains. 2-rank morphism means the rank of
source object and target object differ by 2. In other words, a
\(1\)-dimensional facet should have exactly 2 \(0\)-dimensional
subfacets, that means only line segments are valid \(1\)-dimensional
facets. Moreover, the multiplication of the sign of chains should be
different, that is, a valid diamond should obey
\(\phi_+ \circ \psi_+ = \phi_- \circ \psi_-\) and
\(|\phi_+| \times |\psi_+| = - |\phi_-| \times |\psi_-|\). This property
captures the topological nature of the Euclidean space. To include
infinite lines and rays as valid \(1\)-dimensional facets, this property
should be loosen. The \textbf{semi-diamond property} can be stated as: a
2-rank morphism can be divided into \emph{at most} two maximum chains.
In other words, a \(1\)-dimensional facet should have at most 2
\(0\)-dimensional subfacets.

Except diamond property, the conventional abstract polytope have two
more constraints: strongly connected and uniform maximal chains.
Strongly connected means all facets in the polytope are firmly glued
together. With this property, it is not valid to glue two cubes together
with one edge. Also, a disk with a hole is also invalid, since its
boundaries are not connected. Uniform maximal chains property states
that all maximal chains contain the same number of facets. With this
property, the ranks of facets of a maximal chain are differed by \(1\)
adjacently, so that a rank function can be derived naturally. That means
all facets only have direct subfacets with dimensions differing by
\(1\). They can be generalized as two additional properties for an
abstract orientable incidence structure. That is, \textbf{strongly
decomposable}: all morphisms (except terminal morphisms) can be
decomposed down to 1-rank morphisms; \textbf{strongly unsplittable}: all
section categories of morphisms with rank greater than 2 are
unsplittable. A weaker version, called \textbf{strongly initial
unsplittable}, can be defined as: a bounded acyclic category is said to
be initial unsplittable iff all section categories of initial morphisms
with rank greater than 2 are unsplittable. Strongly initial unsplittable
further states that section categories of all non-initial morphisms can
always be splitted into strongly initial unsplittable categories.
Strongly initial unsplittable relaxes the constraint of shared vertices,
but still limit to connected boundaries.

If only facets homeomorphic to \(n\)-disks are considered, it becomes a
semi-regular and normal CW complex. A \textbf{regular CW complex} is a
CW complex whose gluing maps are homeomorphisms onto images. A CW
complex is said to be \textbf{normal} if each closed cell is a
subcomplex. The incidence structure of a regular and normal CW complex
is given by inclusion relations between closed cells, which forms a CW
poset \autocite{bjorner1984}. Define \textbf{semi-regular CW complex} as
a CW complex whose gluing maps are local homeomorphisms onto images. A
semi-regular and normal CW complex possesses an incidence structure
manifested by the fibers of gluing maps, which forms a bounded acyclic
category. It obeys diamond property and is strongly decomposable and
strongly initial unsplittable: it is initial unsplittable due to
connectivity of cells, and because of local embedding gluing maps, the
vertex figure of a cell is always a disjoint union of simply-connected
components, which leads to strongly initial unsplittable property. It
seems there are more constraints on the incidence structure for a
semi-regular and normal CW complex, such that the geometric realization
of an abstract orientable incidence structure satisfying such
constraints is also a semi-regular and normal CW complex.

\begin{figure}
\hypertarget{fig:chains}{%
\centering
\includegraphics[width=0.6\textwidth,height=\textheight]{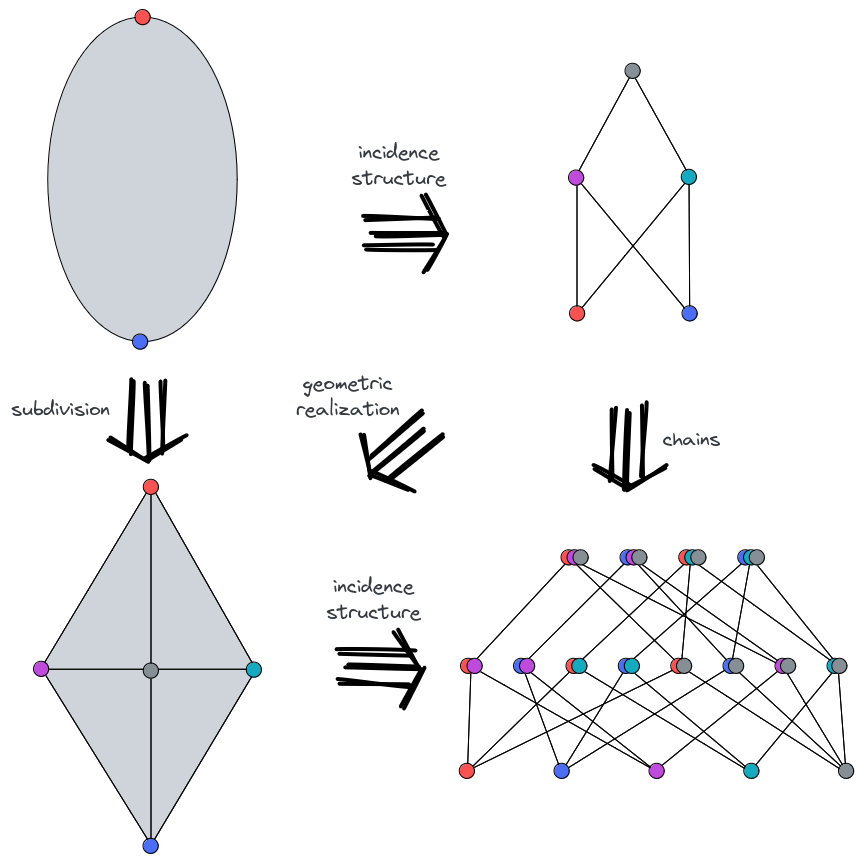}
\caption{Relations between incidence structure and chains. The incidence
structure of a polytope is its face lattice (top left to top right),
which is drawn as a Hasse diagram. The barycentric subdivision of a
polytope is built by dividing facets according to the barycenters (top
left to bottom left), which results in a simplicial complex. Taking
geometric realization of the face lattice (the order complex of a poset)
also results in the same complex (top right to bottom left). The
incidence structure of this complex (bottom left to bottom right) is
just the poset of chains of the face lattice (top right to bottom
right).}\label{fig:chains}
}
\end{figure}

\begin{figure}
\hypertarget{fig:nerve}{%
\centering
\includegraphics[width=0.6\textwidth,height=\textheight]{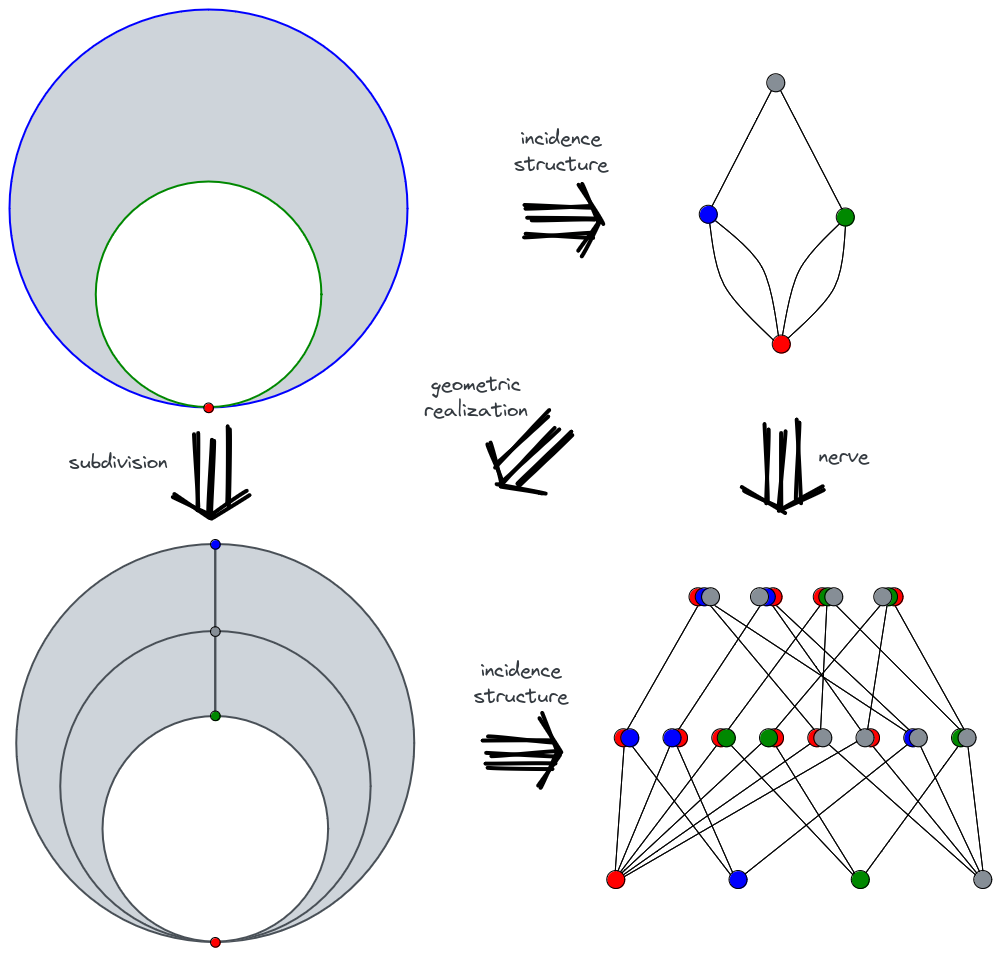}
\caption{Relation between incidence structure and nerve. The incidence
structure of a semi-regular normal CW complex is the bounded acyclic
category formed by fibers of gluing maps (top left to top right). The
barycentric subdivision is defined in a similar way (top left to bottom
left), which results in a geometrical simplicial set. Taking geometric
realization of the incidence structure (the geometric realization of its
nerve) also results in the same complex (top right to bottom left). The
incidence structure of this complex (bottom left to bottom right) is
just the nerve of the bounded acyclic category (top right to bottom
right).}\label{fig:nerve}
}
\end{figure}

Under this constraint, the barycentric subdivision can be defined
without geometric property. Recall that the barycentric subdivision of a
conventional polytope is just the geometric realization (more precisely,
the order complex) of its face lattice (see Figure \ref{fig:chains}).
The barycentric subdivision of a semi-regular normal CW complex can be
defined similarly, and it is just the geometric realization of the nerve
of its incidence structure. In our definition of morphism chains, it
should be defined as: the geometric realization of a nerve is
constructed by making \((n{-}1)\)-simplex for each non-degenerated
\(n\)-chain \(\langle \phi_n, \dots, \phi_2, \phi_1 \rangle\) which is
bounded by the corresponding simplexes of its direct subchains, where
the orientation of this simplex should be defined such that the first
boundary \(\langle \phi_n \circ \phi_{n-1}, \dots \rangle\) is
positively-oriented with respect to it, and the second boundary
\(\langle \phi_n, \phi_{n-1} \circ \phi_{n-2}, \dots \rangle\) is
negatively-oriented with respect to it, and so on. The \((-1)\)-simplex
is the null face. The geometric realization of the nerve of a
semi-regular normal CW complex is just its barycentric subdivision (see
Figure \ref{fig:nerve}). More interestingly, the incidence structure of
a vertex figure is just making a vertex figure of the geometric
realization of upper closure of its incidence structure. There is a
similar relation to the face figure. It shows \textbf{``the incidence
structure of a geometric object is related to the geometry of an
incidence structure''}. This series of articles will not discuss such
correspondence further, since it only appears in a restricted incidence
structure. In the next article, we will focus on finite bounded acyclic
categories, and develop algorithms for them.

\printbibliography

\end{document}